\documentclass[11pt,twoside]{article}
\usepackage{times}
\usepackage{amsmath,amssymb,amsthm}
\usepackage{enumerate}
\usepackage{cite}
\usepackage{graphicx}
\usepackage{booktabs}
\usepackage{makecell}
\usepackage{epstopdf}
\usepackage{xcolor}
\usepackage[hidelinks]{hyperref}

\allowdisplaybreaks

\newcommand{\supp}{\operatorname{supp}}
\newcommand{\R}{\mathbb R}
\newcommand{\p}{\partial}
\newcommand{\ve}{\varepsilon}
\newcommand{\f}{\frac}
\newcommand{\na}{\nabla}
\newcommand{\al}{\alpha}
\newcommand{\Ecal}{\mathcal E}

\newcommand{\dd}{\,\mathrm{d}}
\newcommand{\crit}{\textup{crit}}

\theoremstyle{plain}
\newtheorem{theorem}{Theorem}[section]
\newtheorem{proposition}{Proposition}[section]
\newtheorem{lemma}[theorem]{Lemma}

\theoremstyle{definition}

\theoremstyle{remark}
\newtheorem{remark}{Remark}[section]
\numberwithin{equation}{section}

\title{Global existence of small data solutions to 3-D semilinear
	Euler-Poisson-Darboux equations, II}

\author{Li Qianqian$^{1}$,\quad Xin Qiao$^{2}$,\quad Yin Huicheng$^{3,}$\footnote{Li Qianqian (\texttt{214597007@qq.com}) and Yin Huicheng
    (\texttt{huicheng@} \texttt{nju.edu.cn}, \texttt{05407@njnu.edu.cn}) are supported by the National Natural Science Foundation of China
    (No.~12331007, No.~12571237).  Xin Qiao (\texttt{xinqiaoylsy@163.com}) is supported by the National Natural Science Foundation of China
    (No.~12261092) and the Scientific
Research Program of the Higher Education Institution of XinJiang (No. XJEDU2025P090).}\vspace{0.5cm}\\
\small
1. College of Mathematics and System Science, Xinjiang University, Urumqi 830017, China.\\
\small
2. School of Mathematics and Statistics, Yili Normal University, Yining, 835000, China.\\
\small
3. School of Mathematical Sciences and IMS, Nanjing Normal University, Nanjing 210023, China.\\
}
\vspace{0.5cm}

\author{Li Qianqian$^{1}$,\quad
	Yin Huicheng$^{2,}$\footnote{Li Qianqian
		(\texttt{214597007@qq.com}) and Yin Huicheng
		(\texttt{huicheng@nju.edu.cn}, \texttt{05407@njnu.edu.cn})
		are supported by the National Natural Science Foundation of China
		(No.~12331007, No.~12571237). }\vspace{0.5cm}\\
	\small
	1. College of Mathematics and System Science, Xinjiang University,
	Urumqi 830017, China.\\
	\small
	2. School of Mathematical Sciences and IMS, Nanjing Normal University,
	Nanjing 210023, China.\\
}

\begin{document}
	\date{}
	\maketitle
	\thispagestyle{empty}
	
\begin{abstract}
For the 3-D semilinear Euler-Poisson-Darboux equation $\square u+\f{\mu}{t}\p_tu=|u|^p$, where $t\geq1$,
$\mu>0$ and $p>1$, it is conjectured that the small data solution $u$ exists globally when $p>p_s(3+\mu)=\f{\mu+4+\sqrt{\mu^2+16\mu+32}}{2(\mu+2)}$
for $0<\mu\le \f{14}{5}$, and when $p>p_f(3)=\f53$ for $\mu\geq\f{14}{5}$, where $p_s$ and $p_f$ denote the Strauss and Fujita
exponents, respectively. In our previous paper [Li Qianqian, Yin Huicheng, Global existence of small data solutions to
3-D semilinear Euler-Poisson-Darboux equations. Preprint (2026), arXiv:2607.04575], for  $p>\max\{\frac{5}{3},1+\frac{2}{\mu}\}$ and $\mu\geq\frac{14}{5}$, the global small solution $u$ was obtained. In this paper, we establish the
global existence of $u$ for the remaining range of $\f{5}{3}<p\leq1+\f{2}{\mu}$ and $\mu\ge \f{14}{5}$.
Therefore, the above 3-D conjecture for $\mu\geq\f{14}{5}$ is solved.
\end{abstract}
	
	
	\vskip 0.1 true cm
	
	
	\tableofcontents
	
	%% -------------------------------------------------------------------
	\section{Introduction}\label{Sec1}
	
	In this paper, we continue to study the 3-D semilinear Euler-Poisson-Darboux equation
	\begin{equation}\label{qq:1}
		\left\{
		\begin{aligned}
			&\square u+\f{\mu}{t}\p_tu=|u|^p,
			&&t\geq1,\quad x\in\R^3,\\
			&(u,\p_tu)(1,x)=\ve(u_0,u_1)(x),
		\end{aligned}
		\right.
	\end{equation}
	where $\mu>0$, $p>1$, $\square=\p_t^2-\Delta$,
	$\Delta=\p_1^2+\p_2^2+\p_3^2$, $\p_i=\p_{x_i}$
	$(i=1,2,3)$, $(u_0,u_1)\in C_0^\infty(\R^3)
	\times C_0^\infty(\R^3)$ with
	$\supp(u_0,u_1)\subset B(0,1)$, and $\ve>0$ is
	small. For some interesting  physical backgrounds for the
	linear operator $\square+\f{\mu}{t}\p_t$ and semilinear wave equation in 2-D or 3-D problem \eqref{qq:1}, one can see
	\cite{LWY},\cite{PR-0} and the references therein.
	
	For $x\in\Bbb R^n$ ($n\geq 2$) and $\mu>0$, set
	\[
	p_{\mathrm{crit}}(n,\mu)=\max\{p_s(n+\mu),p_f(n)\},
	\]
	where
	$
	p_s(z)=\frac{z+1+\sqrt{z^2+10z-7}}{2(z-1)}
	\, (z>1)
	$
	stands for the Strauss exponent, which is the positive root of the quadratic algebraic equation
	$(z-1)p^2-(z+1)p-2=0$,
	and $p_f(n)=1+\frac{2}{n}$ stands for the Fujita exponent.
	Note that the Strauss exponent $p_s(n)$ originates from the critical exponent
	problem for the global existence/blowup of solution $v$ to the
	$n$-dimensional semilinear wave equation
	$\square v=|v|^p$ (see \cite{Strauss}), while the Fujita exponent
	$p_f(n)$  stems from the corresponding $n$-dimensional semilinear parabolic
	equation $\partial_t v-\Delta v=|v|^p$ (see \cite{Fuj}).
	In particular, it holds that
	\[
	p_{\mathrm{crit}}(3,\mu)=
	\begin{cases}
		p_s(3+\mu), & 0<\mu\leq \frac{14}{5},\\[2mm]
		p_f(3)=\frac{5}{3}, & \mu\geq\frac{14}{5}.
	\end{cases}
	\]
	In view of \cite{Rei1,Imai} and \cite{LY-2}, with respect to the 3-D problem \eqref{qq:1} with more physical significance and background,
there is an open question as follows:
	
	\vskip 0.1 true cm
	\noindent
	{\bf Open question (A).}
	{\it For the problem \eqref{qq:1},
		
		\noindent
		{\bf (A1)} if $\mu\geq\f{14}{5}$ and
		$p>p_{\crit}(3,\mu)=\f53$, there is a global small solution $u$;
		
		\noindent
		{\bf (A2)} if $0<\mu<\f{14}{5}$ and
		$p>p_{\crit}(3,\mu)=p_s(3+\mu)$, the global small solution $u$ exists.
	}

	We now recall some interesting results concerning the
	critical exponent problem for the
	corresponding $n$-dimensional problem \eqref{qq:1}.
	The blowup result has been completely solved  for
	$1<p\leq p_{\crit}(n,\mu)$ (see
	\cite{FLT,IS,LTW,PR-0,PR,TL2,W1}).
	However, so far no systematic results on global existence have been obtained for
$p>p_{\crit}(n,\mu)$ and $n\ge 3$. A summary of the relevant work can be found
	in \cite{LY-2}, and further references are given in
	\cite{DA-0,DA,Rei2,Rei1,HL,HLY,LZ-1,LWY,LY,P1}.
	
	In our previous paper \cite{LY-2}, we have established the
	global existence of a small solution to 3-D problem \eqref{qq:1}
	when
	\begin{equation}\label{qq:2}
		\text{$p>\max\big\{\f53,1+\f2\mu\big\}$
		and $\mu\geq\f{14}{5}$}.
	\end{equation}
	Due to
	\begin{equation*}\label{qq:3}
		\max\big\{\f53,1+\f2\mu\big\}=\f53
		\quad\text{for }\mu\geq3,
	\end{equation*}
	Theorem~1.1 of \cite{LY-2} obviously yields the global existence
	for {\bf Open question (A1)} with $\mu\geq3$.
	However, for $\f{14}{5}\leq\mu<3$, the condition
	\eqref{qq:2} requires the restriction of $p>1+\f2\mu$. This means that
{\bf Open question (A1)} is not solved for the remaining range of
	\begin{equation}\label{qq:4}
		\f53<p\leq1+\f2\mu\quad\text{and $\f{14}{5}\leq\mu<3$}.
	\end{equation}
	It is pointed out that in \cite{LQY}, under the condition \eqref{qq:4} and
the radial solution case of \eqref{qq:1}, the global existence
	of small solution has been established by making use of the hypergeometric Riemann
	representation and deriving some new weighted spacetime pointwise
	estimates. In order to solve {\bf Open question (A1)} completely, it is required to remove
	the radial symmetry condition in \cite{LQY}.
	In this paper, based on the positivity of fundamental solution to 3-D linear wave equation, by constructing suitable radial
	auxiliary solutions and developing a comparison
	argument to control the corresponding Picard iteration for problem \eqref{qq:1},
	we can complete the proof on the global existence of small solution $u$ to \eqref{qq:1}
	under \eqref{qq:4}. Our main result is stated as follows.
	
	\begin{theorem}\label{thm:main}
		For $\f53<p\leq1+\f2\mu$ and
		$\f{14}{5}\leq\mu<3$, there exists a constant
		$\ve_0>0$ such that when $0<\ve\leq\ve_0$,
		problem \eqref{qq:1} admits a unique global solution
		\begin{equation}\label{qq:5}
			u\in C([1,\infty);H^2(\R^3))
			\cap C^1([1,\infty);H^1(\R^3))
			\cap C^2([1,\infty);L^2(\R^3)).
		\end{equation}
		Meanwhile, the solution $u$ satisfies the following
		decay estimate:
		\begin{equation}\label{qq:6}
			\|u(t,\cdot)\|_\infty
			\leq C\ve\,t^{-\frac{\mu}{2}}(1+t)^{-1},
			\quad t\geq1,
		\end{equation}
		where $C>0$ is a generic constant depending on $\mu, p, u_0$ and $u_1$.
	\end{theorem}
	
	\begin{remark}\label{re1}
		Combining Theorem~\ref{thm:main} with the global
		existence result in Theorem~1.1 of \cite{LY-2} for
		$p>\max\{\f53,1+\f2\mu\}$ with
		$\mu\geq\f{14}{5}$, {\bf Open question (A1)} is
		systematically solved.
	\end{remark}
	
	\begin{remark}\label{re-2}
		In Theorem~\ref{thm:main}, we have removed the radial symmetry
		assumption on the initial data in \cite{LQY}. However,
the result in \cite{LQY} still plays a basic role in proving Theorem~\ref{thm:main}
(see Proposition \ref{prop:radial} and the choice of initial data in \eqref{qq:158} for finding
a suitable comparison function below).
	\end{remark}
	
	\begin{remark}\label{re-3}
		It is pointed out that the global existence of small
		data solutions to problem \eqref{qq:1} has been
		established for $\mu\in(0,1)\cup(1,2)$ and
		$p>p_{\crit}(3,\mu)$, see
		\cite[Remark~1.14]{HLY}. Therefore,
		{\bf Open question (A2)} only remains open for  $p>p_{\crit}(3,\mu)$ with
		$\mu\in\{1\}\cup (2,\f{14}{5})$. This will be studied in our forthcoming paper by
applying the Bessel functions and establishing some crucial Strichartz estimates as in \cite{LWY}.
When this work is finished, then the {\bf Open question (A2)} will be solved completely.
	\end{remark}

\begin{remark}\label{re-4}
From the decay estimate \eqref{qq:6}, it is to see that this decay rate is much better than
the corresponding rates for solutions of 3D wave equation (decay rate $(1+t)^{-1}$), the 3D Klein-Gordon  equation (decay rate $(1+t)^{-3/2}$),
or 3D parabolic equation (decay rate $(1+t)^{-3/2}$).  Similarly, as explained for $\mu=2$ in problem \eqref{qq:1},
the decay behavior in \eqref{qq:6} corresponds exactly to the same
decay as that of the wave equation in dimension ``$3+\mu$".
\end{remark}

\begin{remark}\label{re-5}
The methods in our series of papers (\cite{HLY} and \cite{LWY}-\cite{LQY}) can be utilized to study the more general problem
\begin{equation}\label{qq:1-1}
		\left\{
		\begin{aligned}
			&\square u+\f{\mu}{t}\p_tu=t^{\alpha}|u|^p,
			&&t\geq1,\quad x\in\R^n,\\
			&(u,\p_tu)(1,x)=\ve(u_0,u_1)(x),
		\end{aligned}
		\right.
	\end{equation}
where $n\ge 2$ and $\alpha\in\Bbb R$. Note that for the regular $n$-dimensional semilinear Euler-Poisson-
Darboux-Tricomi equation problem
\begin{equation}\label{qq:1-2}
		\left\{
		\begin{aligned}
			&\p_t^2 u-t^{2m}\Delta u+\f{\mu}{t}\p_tu+\f{\nu}{t^2}u=|u|^p,
			&&t\geq1,\quad x\in\R^n,\\
			&(u,\p_tu)(1,x)=\ve(u_0,u_1)(x),
		\end{aligned}
		\right.
	\end{equation}
where $m>-1$, $n\ge 2$, $\mu>0$ and $\nu>0$, under some suitable transformations of the unknown function $u$
and the time variable $t$, the problem \eqref{qq:1-2} can be transformed into problem \eqref{qq:1-1}
(see \cite{PR-0}).

\end{remark}

	We now sketch the proof of Theorem~\ref{thm:main}.
	Following the physical space approach in
	\cite{Asakura1986,Rei1}, we can remove the
	first order time derivative term $\f{\mu}{t}\p_tu$
in \eqref{qq:1} by the Liouville
	transformation $v(t,x)=t^{\frac{\mu}{2}}u(t,x)$.
	In this case, problem \eqref{qq:1} admits the following form
	\begin{equation}\label{qq:7}
		\left\{
		\begin{aligned}
			&\square v=\f{\mu(\mu-2)}{4t^2}v
			+t^{-\frac{\mu(p-1)}{2}}|v|^p,
			&&t\geq1,\quad x\in\R^3,\\
			&(v,\p_tv)(1,x)=\ve(f,g)(x),
		\end{aligned}
		\right.
	\end{equation}
	where $f=u_0$ and $g=u_1+\f\mu2u_0$.
	Note that for the nonradial function $v=v(t,x)=v(t,r,\omega)$, it holds
	\begin{equation}\label{qq:8}
		r\square v=(\p_t^2-\p_r^2)(rv)
		-\f1{r^2}\Delta_{\mathbb S^2}(rv),\quad\text{where $r=|x|>0$ and $\omega=\f{x}{r}\in\mathbb S^2$}.
	\end{equation}
Due to  $\Delta_{\mathbb S^2}(rv)=0$ for radial $v=v(t,r)$,
then for the symmetric solution of problem \eqref{qq:7},
it can be transformed into the corresponding 1-D semilinear equation $(\p_t^2-\p_r^2)(rv)=\f{\mu(\mu-2)}{4t^2}rv
			+t^{-\frac{\mu(p-1)}{2}}r^{1-p}|rv|^p$ of $rv$.
	However, for the non-radial solution $v$, the angular term $\f1{r^2}\Delta_{\mathbb S^2}v$
naturally appears in $\Delta v$. To overcome this basic difficulty, we will
	construct a suitable radial auxiliary solution to
	control the non-radial Picard iteration for problem \eqref{qq:7}.
	The main ingredients are summarized as follows.
	
	\vskip 0.1 true cm
	
	${\bf \bullet}$ By the Kirchhoff
	formula for 3-D linear wave equation $\square z=h(t,x)$ with
$(z,\p_tz)(1,x)=(f,g)(x)$, one has
$$z=H(f,g)+I(h),$$
where
\begin{equation}\label{qq:9}
\begin{aligned}
&H(f,g)=\f1{4\pi}\int_{\mathbb S^2}
\bigl[f(x+(t-1)\omega)
+(t-1)\na f(x+(t-1)\omega)\cdot\omega+(t-1)g(x+(t-1)\omega)\bigr]\dd\omega,\\
&I(h)=\f1{4\pi}\int_1^t(t-s)
\int_{\mathbb S^2}h(s,x+(t-s)\omega)\dd\omega\dd s.
\end{aligned}
\end{equation}
Hence, $I(h)\ge 0$ holds when $h\ge 0$.
This observation will be applied in the comparison argument below.
	
	\vskip 0.1 true cm
	
	${\bf \bullet}$ For $\f53<q\leq1+\f2\mu$, it follows from
	the radial result in \cite{LQY} that $U$ solves the following problem globally
	\begin{equation}\label{qq:12}
		\left\{
		\begin{aligned}
			&\square U+\f\mu t\p_tU=|U|^q,
			\qquad t\geq1,\quad x\in\R^3,\\
			&(U,\p_tU)(1,x)=\ve\big(F,G-\f\mu2F\big)(x),
		\end{aligned}
		\right.
	\end{equation}
	where
	$F,G\in C_0^\infty(B(0,1))$ are some suitable radial functions to be chosen later.
	Moreover, one has that $W=t^{\frac{\mu}{2}}U$ satisfies
	\begin{equation}\label{qq:13}
		W=\ve H(F,G)
		+I\big(\frac{\mu(\mu-2)}{4s^2}W
		+s^{-\frac{\mu(q-1)}{2}}|W|^q\big)\quad\text{and \quad$|W|\leq C\ve(1+t)^{-1}$}.
	\end{equation}
	Therefore, when $H(F,G)\geq0$, it also holds
	\begin{equation}\label{qq:14}
	W(t,x)\ge 0.
	\end{equation}
	
	\vskip 0.1 true cm
	
	${\bf \bullet}$ Choosing the nonnegative radially symmetric functions
	$(F,G)$ in \eqref{qq:12} as
	\begin{equation}\label{qq:15}
		F\geq|f|,\quad
		G\geq|g|+|\na f|+|\na F|.
	\end{equation}
Then it follows from the expression in \eqref{qq:9} that
\begin{equation}\label{qq:16}
	H(F,G)(t,x)\geq|H(f,g)(t,x)|.
\end{equation}

	\vskip 0.1 true cm
	
	${\bf \bullet}$ Let the number $q$ satisfy
	\begin{equation}\label{qq:18}
	\f53<q<\min\big\{p,1+\f2\mu\big\}.
	\end{equation}
	For small $\ve>0$, one has $0\leq s^{-\frac{\mu}{2}}W(s,x)\leq1$ by \eqref{qq:13}-\eqref{qq:14}.
	This, together with $q<p$, yields
	\begin{equation}\label{qq:19}
		\begin{aligned}
			I\big(\f{\mu(\mu-2)}{4s^2}W
			+s^{-\frac{\mu(p-1)}{2}}W^p\big)&=I\big(\f{\mu(\mu-2)}{4s^2}W
			+s^{-\frac{\mu(q-1)}{2}}W^q
			\bigl(s^{-\frac{\mu}{2}}W\bigr)^{p-q}\big)\\&\leq
			I\big(\f{\mu(\mu-2)}{4s^2}W
			+s^{-\frac{\mu(q-1)}{2}}W^q\big).
		\end{aligned}
	\end{equation}
	Then it follows from \eqref{qq:16}  and \eqref{qq:13} that
	\begin{equation}\label{qq:20}
		|\ve H(f,g)|
		+I\big(\f{\mu(\mu-2)}{4s^2}W
		+s^{-\frac{\mu(p-1)}{2}}W^p\big)\leq W.
	\end{equation}
In order to solve \eqref{qq:7}, we take the following Picard iteration scheme $(k\geq0)$
	\begin{equation}\label{qq:17}
			\begin{aligned}
			&v_{k+1}=\ve H(f,g)
			+I\big(\f{\mu(\mu-2)}{4s^2}v_k
			+s^{-\frac{\mu(p-1)}{2}}|v_k|^p\big),
			\end{aligned}
		\end{equation}
where $v_0=0$. By \eqref{qq:20} and induction method, we can derive the uniform decay estimate
	from \eqref{qq:17}
	\begin{equation}\label{qq:21}
		|v_k(t,x)|\leq W(t,x)\leq C\ve(1+t)^{-1},\quad k\geq0.
	\end{equation}
	Meanwhile, it can be shown that $v_k\to v$ and $v$ is a solution of \eqref{qq:7}.
	By $u=t^{-\frac{\mu}{2}}v$, then Theorem \ref{thm:main} is derived.

	\vskip 0.2 true cm
	\noindent\textbf{Notations}
	\vskip 0.1 true cm
	
	$\bullet$ For $f,g\geq0$,  denote $f\lesssim g$ when
	$f\leq Cg$, and $f\sim g$ when
	$f\lesssim g$ and $g\lesssim f$,
	where $C>0$ is a generic constant independent of $\ve$.
	
	\vskip 0.1 true cm
	
	$\bullet$ For $x\in\R$, set
	$\langle x\rangle=(1+|x|^2)^{\frac{1}{2}}$ and $|x|=\sqrt{x_1^2+x_2^2+x_3^2}$.
	
	\vskip 0.1 true cm
	
	$\bullet$ $W_-=\max\{-W,0\}.$
	
	\vskip 0.1 true cm
	
	$\bullet$ For $\frac{5}{3}<p\leq1+\frac{2}{\mu}$
	and $\frac{14}{5}\leq\mu<3$, define
	\begin{equation}\label{qq:22}
		\begin{aligned}
			&a=\frac{\mu}{2}\in
			\big[\frac{7}{5},\frac{3}{2}\big),
			\quad
			c=a(a-1)\in
			\big[\frac{14}{25},\frac{3}{4}\big),\\
			&b_m(t)=t^{-a(m-1)}\in(0,1]
			\quad\text{for $t\geq1$ and $m\in\{p,q\}$},
		\end{aligned}
	\end{equation}
	where $q$ is defined in \eqref{qq:18}.

\section{Preliminaries}\label{Sec2}

Applying the Liouville transformation
$v=t^au$, we have that from \eqref{qq:1}
\begin{equation}\label{qq:25}
\left\{
\begin{aligned}
&\square v-\frac{c}{t^2}v=b_p(t)|v|^p,
&&t\geq1,\quad x\in\R^3,\\
&(v,\p_tv)(1,x)=\ve(f,g)(x),
\end{aligned}
\right.
\end{equation}
where $f=u_0$, $g=u_1+au_0$, and $a,c,b_p$
are defined in \eqref{qq:22}.

By the Kirchhoff formula
(see \cite[Section~2.4]{Evans}), the solution $z$
of $\square z=0$ with
$(z,\p_tz)(1,x)=(f,g)(x)$ is given by
\begin{equation}\label{qq:30}
\begin{aligned}
z(t,x)=H(f,g)(t,x)
=&\frac{1}{4\pi}\int_{\mathbb S^2}
\bigl[f(x+(t-1)\omega)
+(t-1)\na f(x+(t-1)\omega)\cdot\omega\\
&\hspace{40mm}+(t-1)g(x+(t-1)\omega)\bigr]\dd\omega.
\end{aligned}
\end{equation}
Meanwhile, for $\square w=h$ with
$(w,\p_tw)(1,x)=(0,0)$, one has
\begin{equation}\label{qq:31}
w(t,x)=I(h)(t,x)=
\frac{1}{4\pi}\int_1^t(t-s)\int_{\mathbb S^2}h(s,x+(t-s)\omega)\dd\omega\dd s.
\end{equation}

From \eqref{qq:31}, one easily has
	
\begin{lemma}\label{lem:positive}
For $T>1$, let $h$ and $\tilde h$ be continuous and bounded
on $[1,T]\times\R^3$.
Then $I(h)\in C([1,T]\times\R^3)$,
\begin{equation}\label{qq:32}
\text{$I(h)\geq0$ for $h\geq0$,
and
$|I(h)|\leq I(|h|)$}.
\end{equation}
Moreover, it holds that
\begin{equation*}\label{qq:33}
\|I(h)(t,\cdot)\|_\infty
\leq\int_1^t(t-s)\|h(s,\cdot)\|_\infty\dd s
\end{equation*}
and
\begin{equation}\label{qq:36}
\tilde h\leq h\Longrightarrow I(\tilde h)\leq I(h).
\end{equation}
\end{lemma}

Assume that
\begin{equation}\label{qq:66}
\frac{14}{5}\leq\mu<3,\quad
\frac{5}{3}<q\leq1+\frac{2}{\mu},\quad
\lambda=a-1,\quad \beta=a(q-1),
\end{equation}
where $a$ has been defined in \eqref{qq:22}.
Choosing $\kappa$ such that
\begin{equation}\label{qq:67}
\frac{2-\beta}{q-1}
<\kappa<\min\{q+\beta-1,\,2-\lambda\}.
\end{equation}
In fact, from \eqref{qq:66}, a direct computation yields
\begin{equation*}\label{qq:69}
2-\lambda-\frac{2-\beta}{q-1}
=\frac{3q-5}{q-1}>0\quad\text{and}\quad
q+\beta-1-\frac{2-\beta}{q-1}>0,
\end{equation*}
which means that the choice of $\kappa$  in \eqref{qq:67} is reasonable.
Note that $q-1\leq\frac{1}{a}$ yields
$\kappa>\frac{2-\beta}{q-1}
=\frac{2}{q-1}-a\geq a>1$. In \cite{LQY}, we have defined the space
$X_\kappa$ which consists of radially symmetric functions
 $w\in C([1,\infty)\times\R)$ satisfying
\[
w(t,-r)=w(t,r),\quad
w(t,r)=0\quad\text{for $|r|>t$},\quad
\p_r(rw)\in C([1,\infty)\times\R).
\]
Meanwhile, the related norm of $X_\kappa$  is
\begin{equation}\label{qq:72}
\|w\|_{X_\kappa}
=\sup_{\substack{t\geq1\\0\leq r\leq t}}
(1+t+r)(1+t-r)^{\kappa-1}
\big(
|w(t,r)|+\frac{|\p_r(rw)(t,r)|}{1+r}
\big)<\infty.
\end{equation}
Note that the continuity of $w$ implies
\begin{equation}\label{qq:73}
\p_r(rw)(t,0)
=\lim_{r\to0}\frac{rw(t,r)}{r}=w(t,0).
\end{equation}

We now cite the following global existence result
for the radial solution case in \cite{LQY}.

\begin{proposition}\label{prop:radial}
Under \eqref{qq:66}--\eqref{qq:67}, let
$\phi_0,\phi_1\in C_0^\infty(B(0,1))$
be radially symmetric.
Then there exists a constant $\delta_*>0$ such that for
$0<\delta\leq\delta_*$, the Euler-Poisson-Darboux problem
\begin{equation}\label{qq:76}
\left\{
\begin{aligned}
&\square U+\frac{\mu}{t}\p_tU=|U|^q,
&t\geq1,\quad x\in\R^3,\\
&(U,\p_tU)(1,x)=\delta(\phi_0,\phi_1)(x)
\end{aligned}
\right.
\end{equation}
admits a global radial solution $U$ such that
$W=t^aU\in X_\kappa$ satisfies
\begin{equation}\label{qq:77}
\begin{aligned}
\|W\|_{X_\kappa}&\leq C\delta,\quad
|W(t,r)|\leq\frac{C\delta}{1+t}.
\end{aligned}
\end{equation}
\end{proposition}

\begin{proof}
The existence of $U$ and the estimate
$\|W\|_{X_\kappa}\leq C\delta$
follow from \cite[Theorem~1.1 and Section~5]{LQY} immediately.
Due to $\kappa>1$, the evenness of $W$ and
\eqref{qq:72} yield
\begin{equation}\label{qq:78}
|W(t,r)|
\leq \frac{\|W\|_{X_\kappa}}
{\langle t+|r|\rangle\langle t-|r|\rangle^{\kappa-1}}
\leq \frac{C\delta}{1+t},
\qquad |r|\leq t.
\end{equation}
In addition, note that $W(t,r)=0$  for $|r|>t$. Thus \eqref{qq:77} holds.
\end{proof}

We next derive the expression of $W$ in
Proposition~\ref{prop:radial} by the
Kirchhoff form of \eqref{qq:30}--\eqref{qq:31}.
For this purpose, we recall the Riemann kernel
$\Ecal_\lambda$ from \cite{P2,LQY}:
\begin{equation}\label{qq:74}
\begin{aligned}
\Ecal_\lambda(t,s,\sigma)
={}&(ts)^{-\lambda}
\bigl((t+s)^2-\sigma^2\bigr)^\lambda\times{}_2F_1\big(
-\lambda,-\lambda;1;
\frac{(t-s)^2-\sigma^2}{(t+s)^2-\sigma^2}
\big),
\end{aligned}
\end{equation}
where $1\leq s\leq t$, $|\sigma|\leq t-s$,  and the Gauss hypergeometric function is expressed as
(see \cite[Section~15.2]{OLBC})
\begin{equation}\label{qq:75}
{}_2F_1(-\lambda,-\lambda;1;z)
=\sum_{m=0}^{\infty}
\frac{(-\lambda)_m^2}{(m!)^2}z^m
\quad\text{for $|z|<1$}
\end{equation}
with $(d)_0=1$ and
$(d)_m=d(d+1)\cdots(d+m-1)$ for $m\geq1$.

\begin{lemma}\label{lem:compatibility}
Let $W$ be defined in Proposition~\ref{prop:radial}.
Then for $t\geq1$ and $x\in\R^3$, it holds that
\begin{equation}\label{qq:79}
\begin{aligned}
W(t,x)
={}&\delta H(\phi_0,\phi_1+a\phi_0)(t,x)+I\big(
\frac{\lambda(\lambda+1)}{s^2}W+b_q(s)|W|^q
\big)(t,x).
\end{aligned}
\end{equation}
\end{lemma}
	
\begin{proof}
Set $\psi(t,r)=rW(t,r)$ for $r\in\R$.
Then it follows from \eqref{qq:76} and $W=t^aU$ that
\begin{equation}\label{qq:80}
\left\{
\begin{aligned}
&(\p_t^2-\p_r^2)\psi
-\frac{\lambda(\lambda+1)}{t^2}\psi
=rb_q(t)|W(t,r)|^q,
&&t\geq1,\quad r\in\R,\\
&(\psi,\p_t\psi)(1,r)
=\delta r(\phi_0,\phi_1+a\phi_0)(|r|),
\end{aligned}
\right.
\end{equation}
where $a,b_q$ are given in \eqref{qq:22}
and $\lambda$ is defined in \eqref{qq:66}.
By \cite[Sections~2 and~5]{LQY}, $\psi$ can be written as
\begin{equation}\label{qq:81}
\psi(t,r)=\psi_{\rm lin}(t,r)
+c_\lambda\int_1^t\int_{r-t+s}^{r+t-s}
\Ecal_\lambda(t,s,r-\rho)\rho b_q(s)|W(s,\rho)|^q
\dd\rho\dd s,
\end{equation}
where $c_\lambda=2^{-(2\lambda+1)}$ and
\begin{equation}\label{qq:82}
\begin{aligned}
\psi_{\rm lin}(t,r)
={}&\frac{\delta}{2}\Big[
(r+t-1)\phi_0(|r+t-1|)
+(r-t+1)\phi_0(|r-t+1|)\Big]\\
&+\delta c_\lambda\int_{r-t+1}^{r+t-1}
\Ecal_\lambda(t,1,r-\rho)
\rho(\phi_1+a\phi_0)(|\rho|)\dd\rho\\
&-\delta c_\lambda\int_{r-t+1}^{r+t-1}
\left.\p_s\Ecal_\lambda(t,s,r-\rho)\right|_{s=1}
\rho\phi_0(|\rho|)\dd\rho.
\end{aligned}
\end{equation}
Since $\phi_0,\phi_1$ are smooth,
\eqref{qq:74} and \eqref{qq:82} imply
$\psi_{\rm lin}\in C^2([1,\infty)\times\R)$.
Note that $rb_q(t)|W(t,r)|^q$ has a continuous
derivative in $r$. Indeed, for $r\ne0$, one has
\[
W(t,r)=\frac{\psi(t,r)}{r},\quad
\p_rW(t,r)=\frac{\p_r\psi(t,r)-W(t,r)}{r}.
\]
Consequently,
\begin{equation}\label{qq:84}
\begin{aligned}
\p_r\bigl(rb_q(t)|W|^q\bigr)
&=b_q(t)\left[|W|^q
+q|W|^{q-2}W(\p_r\psi-W)\right].
\end{aligned}
\end{equation}
In addition, \eqref{qq:72} and \eqref{qq:73} imply
that $W,\p_r\psi$ are continuous and
$\p_r\psi(t,0)=W(t,0)$. Therefore,
\begin{equation}\label{qq:85}
\begin{aligned}
\left.\p_r\bigl(rb_q(t)|W(t,r)|^q\bigr)\right|_{r=0}
&=\lim_{r\to0}
\frac{rb_q(t)|W(t,r)|^q-0}{r}
=b_q(t)|W(t,0)|^q.
\end{aligned}
\end{equation}
Due to $q>1$, \eqref{qq:84} and \eqref{qq:85} yield
$\p_r\bigl(rb_q(t)|W(t,r)|^q\bigr)
\in C([1,\infty)\times\R)$. Note that a direct computation in \eqref{qq:81} gives
$\psi\in C^2([1,\infty)\times\R)$.
Hence, by  the 1-D d'Alembert's formula and Duhamel's principle
(see \cite[Sections~2.4.1 and~2.4.2]{Evans}),
we obtain from \eqref{qq:80} that
\begin{equation}\label{qq:102}
\begin{aligned}
\psi(t,r)
={}&\frac{\delta}{2}\big[
(r+t-1)\phi_0(|r+t-1|)
+(r-t+1)\phi_0(|r-t+1|)\big]\\
&+\frac{\delta}{2}\int_{r-t+1}^{r+t-1}
\rho(\phi_1+a\phi_0)(|\rho|)\dd\rho\\
&+\frac{1}{2}\int_1^t\int_{r-t+s}^{r+t-s}
\rho\big[
\frac{\lambda(\lambda+1)}{s^2}W(s,\rho)
+b_q(s)|W(s,\rho)|^q
\big]\dd\rho\dd s.
\end{aligned}
\end{equation}
On the other hand, for a continuous radial function
$\phi(x)=\phi(r)$ with $r=|x|>0$, it follows from
direct calculation that for $h\geq0$,
\begin{equation}\label{qq:104}
\begin{aligned}
\frac{h}{4\pi}
\int_{\mathbb S^2}\phi(x+h\omega)\dd\omega
&=\frac{h}{4\pi}\int_0^{2\pi}\int_0^\pi
\phi\big(\sqrt{r^2+h^2+2rh\cos\alpha}\big)
\sin\alpha\dd\alpha\dd\varphi\\
&\overset{\theta=\cos\alpha}{=}
\frac{h}{2}\int_{-1}^1
\phi\big(\sqrt{r^2+h^2+2rh\theta}\big)\dd\theta\\
&\overset{
\rho=\sqrt{r^2+h^2+2rh\theta}
}{=}
\frac{1}{2r}\int_{|r-h|}^{r+h}\rho\phi(\rho)\dd\rho\\
&=\frac{1}{2r}\int_{r-h}^{r+h}\rho\phi(|\rho|)\dd\rho,
\end{aligned}
\end{equation}
where the oddness of $\rho\phi(|\rho|)$ is used in the last equality.

For $x\ne0$,  by \eqref{qq:104} it holds that
\begin{equation}\label{qq:104a}
\begin{aligned}
&\frac{\delta}{2r}\big[
(r+t-1)\phi_0(|r+t-1|)
+(r-t+1)\phi_0(|r-t+1|)\big]\\
&=\frac{\delta}{2r}\p_h
\big(\int_{r-h}^{r+h}\rho\phi_0(|\rho|)\dd\rho\big)
\big|_{h=t-1}=\delta\p_h
\big(\frac{h}{4\pi}
\int_{\mathbb S^2}\phi_0(x+h\omega)\dd\omega\big)\big|_{h=t-1}.
\end{aligned}
\end{equation}
Taking $\phi=\phi_1+a\phi_0$ and $h=t-1$
in \eqref{qq:104}, one has
\begin{equation}\label{qq:104b}
\frac{\delta}{2r}\int_{r-t+1}^{r+t-1}
\rho(\phi_1+a\phi_0)(|\rho|)\dd\rho
=\f{\delta(t-1)}{4\pi}\int_{\mathbb S^2}(\phi_1+a\phi_0)(x+(t-1)\omega)\dd\omega.
\end{equation}
Similarly,
\begin{equation}\label{qq:104c}
\begin{aligned}
&\frac{1}{2r}\int_1^t\int_{r-t+s}^{r+t-s}
\rho\big[
\frac{\lambda(\lambda+1)}{s^2}W(s,\rho)
+b_q(s)|W(s,\rho)|^q
\big]\dd\rho\dd s\\
&\quad=\int_1^t\f{(t-s)}{4\pi}\int_{\mathbb S^2}
\big(\frac{\lambda(\lambda+1)}{s^2}W(s,\cdot)
+b_q(s)|W(s,\cdot)|^q\big)(x+(t-s)\omega)\dd\omega\dd s.
\end{aligned}
\end{equation}
Collecting \eqref{qq:30}-\eqref{qq:31}, \eqref{qq:104a}-\eqref{qq:104c}
and \eqref{qq:102} yields
\begin{equation}\label{qq:103}
\begin{aligned}
W(t,x)=\frac{\psi(t,r)}{r}
={}&\delta H(\phi_0,\phi_1+a\phi_0)(t,x)
+I\big(\frac{\lambda(\lambda+1)}{s^2}W+b_q(s)|W|^q\big)(t,x).
\end{aligned}
\end{equation}
For $x=0$, \eqref{qq:22}, \eqref{qq:72}
and \eqref{qq:77} imply
\[
W,\ \frac{\lambda(\lambda+1)}{t^2}W+b_q(t)|W|^q
\in C([1,\infty)\times\R^3)
\cap L^\infty([1,\infty)\times\R^3).
\]
Then it follows from \eqref{qq:30} and Lemma~\ref{lem:positive} that
\[
\begin{aligned}
W(t,0)
&=\lim_{x\to0}\big[
\delta H(\phi_0,\phi_1+a\phi_0)(t,x)
+I\big(\frac{\lambda(\lambda+1)}{s^2}W+b_q(s)|W|^q\big)(t,x)
\big]\\
&=\delta H(\phi_0,\phi_1+a\phi_0)(t,0)
+I\big(\frac{\lambda(\lambda+1)}{s^2}W+b_q(s)|W|^q\big)(t,0).
\end{aligned}
\]
This, together with \eqref{qq:103}, proves \eqref{qq:79}.
\end{proof}
	
Based on Lemma~\ref{lem:compatibility}, we now derive
the following nonnegativity result for the radial
auxiliary solution $W$.

\begin{lemma}\label{lem:radial-nonnegative}
Let $W$ be defined in Lemma~\ref{lem:compatibility}.
If $H(\phi_0,\phi_1+a\phi_0)\geq0$,
then for $t\geq1$ and $x\in\R^3$, it holds that
\begin{equation}\label{qq:106}
W(t,x)\geq0.
\end{equation}
\end{lemma}

\begin{proof}
For any fixed $T>1$, it follows from \eqref{qq:79}, \eqref{qq:32}
and \eqref{qq:66} that for $t\in[1,T]$,
\begin{equation}\label{qq:107}
\begin{aligned}
&W+I\big(\frac{\lambda(\lambda+1)}{s^2}W_-\big)=\delta H(\phi_0,\phi_1+a\phi_0)
+I\big(\frac{\lambda(\lambda+1)}{s^2}(W+W_-)
+b_q(s)|W|^q
\big)\geq0.
\end{aligned}
\end{equation}
In addition, by \eqref{qq:107} and \eqref{qq:32}, one has
\begin{equation}\label{qq:108}
0\leq W_-
\leq I\big(\frac{\lambda(\lambda+1)}{s^2}W_-\big).
\end{equation}
Combining \eqref{qq:108} with \eqref{qq:33} yields
\begin{equation}\label{qq:109}
\begin{aligned}
\|W_-(t,\cdot)\|_\infty
&\leq\big\|
I\big(\frac{\lambda(\lambda+1)}{s^2}W_-\big)
(t,\cdot)\big\|_\infty\leq\lambda(\lambda+1)\int_1^t
(t-s)\|W_-(s,\cdot)\|_\infty\dd s.
\end{aligned}
\end{equation}
Moreover, we obtain from \eqref{qq:77} that
\[
0\leq\|W_-(t,\cdot)\|_\infty
\leq\|W(t,\cdot)\|_\infty
\leq\frac{C\delta}{1+t}.
\]
For $s\in[1,T]$, by the support condition in
\eqref{qq:72} and the  continuity of $W$, one has
\begin{equation}\label{qq:110}
\left|
\|W_-(t,\cdot)\|_\infty-\|W_-(s,\cdot)\|_\infty
\right|
\leq\sup_{|x|\leq T}|W_-(t,x)-W_-(s,x)|
\longrightarrow0
\quad\text{as }s\to t.
\end{equation}
On the other hand, it follows from \eqref{qq:109} and Gronwall's inequality that
$\|W_-(t,\cdot)\|_\infty=0$.
Therefore, \eqref{qq:106} is shown due to the arbitrary choice of $T>1$.
\end{proof}
	
\section{Comparison principle and Picard iteration}
\label{Sec4}

In this section, we first choose nonnegative radial
functions $F$ and $G$ such that $H(F,G)\geq |H(f,g)|$.

\begin{lemma}\label{lem:homogeneous-majorant}
For any $f,g\in C_0^\infty(B(0,1))$. Then there exist nonnegative radially
symmetric functions $F,G\in C_0^\infty(B(0,1))$ such that
\begin{equation}\label{qq:111}
F\geq|f|,\quad
G\geq|g|+|\na f|+|\na F|.
\end{equation}
Moreover, for $t\geq1$ and $x\in\R^3$, one has
\begin{equation}\label{qq:112}
H(F,G)(t,x)\geq|H(f,g)(t,x)|.
\end{equation}
\end{lemma}

\begin{proof}
Due to $f,g\in C_0^\infty(B(0,1))$, we may choose some constants
$r_0,r_1,r_2$ such that
\begin{equation*}\label{qq:113}
0<r_0<r_1<r_2<1,\quad
\supp f\cup\supp g\subset B(0,r_0).
\end{equation*}
Let $\chi\in C^\infty(\R)$ satisfy $0\leq\chi\leq1$ and
$\chi(\tau)=
\begin{cases}
1, & \tau\leq0,\\
0, & \tau\geq1.
\end{cases}$
\quad In addition, set
\begin{equation*}\label{qq:115}
\begin{aligned}
F(x)&=\|f\|_\infty
\chi\big(\frac{|x|^2-r_0^2}{r_1^2-r_0^2}\big),\\
G(x)&=\big(\|g\|_\infty+\|\na f\|_\infty
+\|\na F\|_\infty\big)
\chi\big(\frac{|x|^2-r_1^2}{r_2^2-r_1^2}\big).
\end{aligned}
\end{equation*}
Then $F,G\in C_0^\infty(B(0,1))$ are nonnegative and
radially symmetric. In addition, it holds that
\begin{equation*}\label{qq:116}
\supp F\subset\overline{B(0,r_1)},\quad
\supp G\subset\overline{B(0,r_2)}.
\end{equation*}
Note that
\begin{equation*}\label{qq:117}
\na F(x)=\frac{2\|f\|_\infty x}{r_1^2-r_0^2}
\chi'\big(\frac{|x|^2-r_0^2}{r_1^2-r_0^2}\big),
\end{equation*}
then
\begin{equation*}\label{qq:118}
\begin{aligned}
|f(x)|
&\leq\|f\|_\infty
\chi\big(\frac{|x|^2-r_0^2}{r_1^2-r_0^2}\big)
=F(x),\\
|g(x)|+|\na f(x)|+|\na F(x)|
&\leq\big(\|g\|_\infty+\|\na f\|_\infty
+\|\na F\|_\infty\big)
\chi\big(\frac{|x|^2-r_1^2}{r_2^2-r_1^2}\big)
=G(x),
\end{aligned}
\end{equation*}
which proves \eqref{qq:111}.

For $h=t-1\geq0$ and $\omega\in\mathbb S^2$,
it follows from \eqref{qq:111} that
\begin{equation}\label{qq:119}
\begin{aligned}
&F(x+h\omega)+h\na F(x+h\omega)\cdot\omega
+hG(x+h\omega)\\
&\geq F(x+h\omega)
+h\bigl(G(x+h\omega)-|\na F(x+h\omega)|\bigr)\\
&\geq |f(x+h\omega)|+h|\na f(x+h\omega)|
+h|g(x+h\omega)|\\
&\geq
\bigl|f(x+h\omega)+h\na f(x+h\omega)\cdot\omega
+hg(x+h\omega)\bigr|.
\end{aligned}
\end{equation}
Integrating \eqref{qq:119} over $\mathbb S^2$
and using \eqref{qq:30}, we arrive at
\begin{equation*}\label{qq:120}
\begin{aligned}
H(F,G)(t,x)
&\geq\frac{1}{4\pi}\int_{\mathbb S^2}
\bigl|f(x+(t-1)\omega)
+(t-1)\na f(x+(t-1)\omega)\cdot\omega
+(t-1)g(x+(t-1)\omega)\bigr|\dd\omega\\
&\geq|H(f,g)(t,x)|.
\end{aligned}
\end{equation*}
which derives \eqref{qq:112}.
\end{proof}
	
Based on Lemma~\ref{lem:homogeneous-majorant}, we next
establish the existence of a global solution
to \eqref{qq:25} by Picard iteration, which will be applied for the proof of Theorem~\ref{thm:main} in Section~\ref{Sec5} below.

\begin{proposition}\label{prop:iteration}
Let $p>1$, $\ve>0$ and $f,g\in C_0^\infty(\R^3)$.
Suppose that a nonnegative function
$W\in C([1,\infty)\times\R^3)$ is bounded on
$[1,T]\times\R^3$ for any $T>1$ and satisfies
\begin{equation}\label{qq:121}
W\geq\ve|H(f,g)|
+I\big(\frac{\lambda(\lambda+1)}{s^2}W+b_p(s)W^p
\big).
\end{equation}
Then \eqref{qq:25} admits a unique global solution
$v$. Moreover, for $t\geq1$ and $x\in\R^3$, it holds that
\begin{equation}\label{qq:122}
|v(t,x)|\leq W(t,x).
\end{equation}
\end{proposition}

\begin{proof}
Define the Picard iteration by
\begin{align}\label{qq:124}
&v_{k+1}=\ve H(f,g)
+I\big(\frac{\lambda(\lambda+1)}{s^2}v_k+b_p(s)|v_k|^p
\big),
\quad k=0,1,2,\ldots,
\end{align}
where $v_0=0$.
It is asserted that for $k\geq0$, $T>1$,
$t\in[1,T]$ and $x\in\R^3$,
\begin{equation}\label{qq:125}
\begin{gathered}
v_k\in C([1,T]\times\R^3)
\cap L^\infty([1,T]\times\R^3),\\
|v_k(t,x)|\leq W(t,x).
\end{gathered}
\end{equation}
We next prove \eqref{qq:125} by induction method.
For $k=0$, one has $v_0=0$ and $|v_0|\leq W$.
Assume that \eqref{qq:125} holds for some $k\geq0$.
By \eqref{qq:22} and the induction hypothesis,
it holds that
\[
\frac{\lambda(\lambda+1)}{s^2}v_k+b_p(s)|v_k|^p
\in C([1,T]\times\R^3)
\cap L^\infty([1,T]\times\R^3).
\]
In view of \eqref{qq:30}, \eqref{qq:124}
and Lemma~\ref{lem:positive}, one then has
$v_{k+1}\in C([1,T]\times\R^3)$.
Moreover, combining \eqref{qq:124} with
\eqref{qq:32} and \eqref{qq:36} derives
\begin{equation}\label{qq:126}
\begin{aligned}
|v_{k+1}|
&\leq\ve|H(f,g)|
+I\big(\frac{\lambda(\lambda+1)}{s^2}|v_k|
+b_p(s)|v_k|^p
\big)\\
&\leq\ve|H(f,g)|
+I\big(
\frac{\lambda(\lambda+1)}{s^2}W+b_p(s)W^p
\big)\\
&
\leq W,
\end{aligned}
\end{equation}
where the last inequality follows from \eqref{qq:121}.
Since $W$ is bounded on $[1,T]\times\R^3$,
we also have $v_{k+1}\in L^\infty([1,T]\times\R^3)$.
Thus the assertion \eqref{qq:125} is shown.

Note that for any $T>1$ and $t\in[1,T]$, and for $k\geq1$,
\begin{equation}\label{qq:130}
v_{k+1}-v_k
=I\big(
\frac{\lambda(\lambda+1)}{s^2}(v_k-v_{k-1})
+b_p(s)(|v_k|^p-|v_{k-1}|^p)
\big).
\end{equation}
Set
\begin{equation}\label{qq:127}
M_T=1+\|W\|_{L^\infty([1,T]\times\R^3)}<\infty.
\end{equation}
Then  it follows from \eqref{qq:130}, \eqref{qq:33} and \eqref{qq:125} that
\begin{equation}\label{qq:131}
\begin{aligned}
&\|v_{k+1}(t,\cdot)-v_k(t,\cdot)\|_\infty\\
&\leq\int_1^t(t-s)
\big\|
\frac{\lambda(\lambda+1)}{s^2}(v_k-v_{k-1})(s,\cdot)
+b_p(s)(|v_k|^p-|v_{k-1}|^p)(s,\cdot)
\big\|_\infty\dd s\\
&\leq\int_1^t(t-s)
\big(\frac{\lambda(\lambda+1)}{s^2}
+pM_T^{p-1}b_p(s)\big)
\|(v_k-v_{k-1})(s,\cdot)\|_\infty\dd s\\
&\leq\bigl(\lambda(\lambda+1)+pM_T^{p-1}\bigr)
\int_1^t(t-s)
\|(v_k-v_{k-1})(s,\cdot)\|_\infty\dd s.
\end{aligned}
\end{equation}
Next we show that for $k\geq0$ and $t\in[1,T]$,
\begin{equation}\label{qq:132}
\|v_{k+1}(t,\cdot)-v_k(t,\cdot)\|_\infty
\leq M_T
\frac{\bigl(\lambda(\lambda+1)+pM_T^{p-1}\bigr)^k
(t-1)^{2k}}{(2k)!}.
\end{equation}
For $k=0$, it follows from
\eqref{qq:124}--\eqref{qq:125} and \eqref{qq:127} that
\begin{equation}\label{qq:133}
\|v_1(t,\cdot)-v_0(t,\cdot)\|_\infty
\leq\|W(t,\cdot)\|_\infty\leq M_T.
\end{equation}
Suppose that \eqref{qq:132} has been proved
for some $k\geq0$.
Next we derive that \eqref{qq:132} also holds for $k+1$.
Combining \eqref{qq:131} with the induction hypothesis
\eqref{qq:132}, one has
\begin{equation}\label{qq:134}
\begin{aligned}
&\|v_{k+2}(t,\cdot)-v_{k+1}(t,\cdot)\|_\infty
\\&\leq
\frac{M_T\bigl(\lambda(\lambda+1)+pM_T^{p-1}\bigr)^{k+1}}
{(2k)!}
\int_1^t(t-s)(s-1)^{2k}\dd s\\
&=
M_T\frac{\bigl(\lambda(\lambda+1)+pM_T^{p-1}\bigr)^{k+1}
(t-1)^{2k+2}}{(2k+2)!}.
\end{aligned}
\end{equation}
Thus the assertion \eqref{qq:132} is shown.
Moreover, a direct computation gives
\begin{equation}\label{qq:138}
\begin{aligned}
\sum_{k=0}^{\infty}
\frac{\bigl(\lambda(\lambda+1)+pM_T^{p-1}\bigr)^k
(T-1)^{2k}}{(2k)!}
&\leq
e^{\bigl(\lambda(\lambda+1)+pM_T^{p-1}\bigr)(T-1)^2}
<\infty.
\end{aligned}
\end{equation}
This, together with
\eqref{qq:132}, yields that for any integer $m>k$,
\begin{equation}\label{qq:139}
\begin{aligned}
\|v_m-v_k\|_{L^\infty([1,T]\times\R^3)}
&\leq
M_T\sum_{j=k}^{\infty}
\frac{\bigl(\lambda(\lambda+1)+pM_T^{p-1}\bigr)^j
(T-1)^{2j}}{(2j)!}
\\
&\longrightarrow0 \qquad \text{as}\, k\to\infty.
\end{aligned}
\end{equation}
Consequently, $v_k$ converges uniformly on
$[1,T]\times\R^3$ to some function
$v\in C([1,T]\times\R^3)
\cap L^\infty([1,T]\times\R^3)$ satisfying
\begin{equation}\label{qq:140}
|v(t,x)|\leq W(t,x).
\end{equation}
In addition, by \eqref{qq:33} and the mean value theorem,
it follows from \eqref{qq:125}, \eqref{qq:127}
and \eqref{qq:140} that
\begin{equation}\label{qq:141}
\begin{aligned}
&\big\|
I\big(
\frac{\lambda(\lambda+1)}{s^2}(v_k-v)
+b_p(s)(|v_k|^p-|v|^p)
\big)
\big\|_{L^\infty([1,T]\times\R^3)}\\
&\leq
\sup_{1\leq t\leq T}\int_1^t(t-s)
\big(
\frac{\lambda(\lambda+1)}{s^2}
+pM_T^{p-1}b_p(s)
\big)
\|v_k(s,\cdot)-v(s,\cdot)\|_\infty\dd s\\
&\leq
\frac{\bigl(\lambda(\lambda+1)+pM_T^{p-1}\bigr)(T-1)^2}{2}
\|v_k-v\|_{L^\infty([1,T]\times\R^3)}
\longrightarrow0  \qquad \text{as}\, k\to\infty.
\end{aligned}
\end{equation}
Passing to the limit in \eqref{qq:124}, we then obtain
from \eqref{qq:141} and the uniform convergence of $v_k$ that
\begin{equation}\label{qq:142}
\begin{aligned}
v
&=\lim_{k\to\infty}v_{k+1}=\ve H(f,g)
+\lim_{k\to\infty}I\big(\frac{\lambda(\lambda+1)}{s^2}v_k+b_p(s)|v_k|^p
\big)\\
&=\ve H(f,g)
+I\big(\frac{\lambda(\lambda+1)}{s^2}v+b_p(s)|v|^p\big).
\end{aligned}
\end{equation}
As the sequence $\{v_k\}$ is independent of $T$ and $T>1$ is arbitrary,
$v$ is a global solution of problem \eqref{qq:25}, and \eqref{qq:140} holds
for $t\geq1$.

Next, we start to show the uniqueness of global
solution $v$ to \eqref{qq:25}.
Let $\widetilde v$ be another global solution
with the same initial data.
Analogously to the treatment of \eqref{qq:131}, we have
\begin{equation*}\label{qq:143}
\begin{aligned}
&\|v(t,\cdot)-\widetilde v(t,\cdot)\|_\infty\\
&=
\big\|
I\big(
\frac{\lambda(\lambda+1)}{s^2}(v-\widetilde v)
+b_p(s)(|v|^p-|\widetilde v|^p)
\big)(t,\cdot)
\big\|_\infty\\
&\leq
\int_1^t(t-s)\bigl[
\frac{\lambda(\lambda+1)}{s^2}
+p\,b_p(s)\big(
1+\|v\|_{L^\infty([1,T]\times\R^3)}
+\|\widetilde v\|_{L^\infty([1,T]\times\R^3)}
\big)^{p-1}
\bigr]\\&\qquad\qquad\qquad\times
\|v(s,\cdot)-\widetilde v(s,\cdot)\|_\infty\dd s\\
&\leq
\big[
\lambda(\lambda+1)
+p\left(
1+\|v\|_{L^\infty([1,T]\times\R^3)}
+\|\widetilde v\|_{L^\infty([1,T]\times\R^3)}
\right)^{p-1}
\big]
\int_1^t(t-s)
\|v(s,\cdot)-\widetilde v(s,\cdot)\|_\infty\dd s,
\end{aligned}
\end{equation*}
which yields
$\|v(t,\cdot)-\widetilde v(t,\cdot)\|_\infty=0$.
Therefore, the uniqueness of the global solution $v$
to \eqref{qq:25} is shown.

To prove that $v$ satisfies the initial conditions in \eqref{qq:25},
we set $t=1$ in \eqref{qq:142}.
Then one has that by \eqref{qq:30} and \eqref{qq:31},
\begin{equation}\label{qq:149}
\begin{aligned}
v(1,x)
&=\ve H(f,g)(1,x)
+I\big(
\frac{\lambda(\lambda+1)}{s^2}v+b_p(s)|v|^p
\big)(1,x)\\
&=\frac{\ve}{4\pi}\int_{\mathbb S^2}f(x)\dd\omega
+\f{1}{4\pi}\int_1^1(1-s)\int_{\Bbb S^2}
\big(
\frac{\lambda(\lambda+1)}{s^2}v(s,\cdot)
+b_p(s)|v(s,\cdot)|^p
\big)(x+(1-s)\omega)\dd\omega\dd s\\
&=\ve f(x).
\end{aligned}
\end{equation}
We now treat the term $\p_t v(1,x)$.
For $0<h\leq T-1$, it follows from \eqref{qq:33},
\eqref{qq:127} and \eqref{qq:140} that
\begin{equation}\label{qq:151}
\begin{aligned}
&\big\|
I\big(
\frac{\lambda(\lambda+1)}{s^2}v+b_p(s)|v|^p
\big)(1+h,\cdot)
\big\|_\infty\\
&\leq
\int_1^{1+h}(1+h-s)
\big\|
\frac{\lambda(\lambda+1)}{s^2}v(s,\cdot)
+b_p(s)|v(s,\cdot)|^p
\big\|_\infty\dd s\\
&\leq
\int_1^{1+h}(1+h-s)
\bigl(\lambda(\lambda+1)M_T+M_T^p\bigr)\dd s\\
&=
\frac{h^2}{2}
\bigl(\lambda(\lambda+1)M_T+M_T^p\bigr).
\end{aligned}
\end{equation}
On the other hand,  it holds that from \eqref{qq:30},
\begin{equation}\label{qq:153}
\begin{aligned}
&\|H(f,g)(1+h,\cdot)-f-hg\|_\infty\\
&=\frac{1}{4\pi}\big\|
\int_{\mathbb S^2}
\big[
f(\cdot+h\omega)-f-h\na f\cdot\omega
+h\bigl(\na f(\cdot+h\omega)-\na f\bigr)\cdot\omega
+h\bigl(g(\cdot+h\omega)-g\bigr)
\big]\dd\omega
\big\|_\infty\\
&\leq\frac{1}{4\pi}\int_{\mathbb S^2}
\bigl[
\|f(\cdot+h\omega)-f-h\na f\cdot\omega\|_\infty
+h\|\na f(\cdot+h\omega)-\na f\|_\infty
+h\|g(\cdot+h\omega)-g\|_\infty
\bigr]\dd\omega\\
&\leq\frac{h^2}{4\pi}\int_{\mathbb S^2}
\bigl[
\int_0^1(1-\theta)
\bigl\|\omega^{\mathsf T}\na^2f(\cdot+\theta h\omega)\omega\bigr\|_\infty
\dd\theta\\
&\qquad\qquad\quad+
\int_0^1
\bigl\|\na^2f(\cdot+\theta h\omega)\omega\bigr\|_\infty
\dd\theta
+\int_0^1
\bigl\|\na g(\cdot+\theta h\omega)\cdot\omega\bigr\|_\infty
\dd\theta
\bigr]\dd\omega\\
&\leq\frac{1}{4\pi}\int_{\mathbb S^2}
\big[
\frac{h^2}{2}\|\na^2f\|_\infty
+h^2\|\na^2f\|_\infty
+h^2\|\na g\|_\infty
\big]\dd\omega\\
&\leq Ch^2
\bigl(\|\na^2f\|_\infty+\|\na g\|_\infty\bigr).
\end{aligned}
\end{equation}
Combining \eqref{qq:142}, \eqref{qq:149} and \eqref{qq:151}-
\eqref{qq:153} yields
\begin{equation*}\label{qq:154}
\begin{aligned}
&\big\|
\frac{v(1+h,\cdot)-v(1,\cdot)}{h}-\ve g
\big\|_\infty\\
&=
\big\|
\frac{\ve}{h}
\bigl(H(f,g)(1+h,\cdot)-f-hg\bigr)
+\frac{1}{h}
I\big(
\frac{\lambda(\lambda+1)}{s^2}v+b_p(s)|v|^p
\big)(1+h,\cdot)
\big\|_\infty\\
&\leq
\frac{\ve}{h}
\|H(f,g)(1+h,\cdot)-f-hg\|_\infty
+\frac{1}{h}\big\|
I\big(
\frac{\lambda(\lambda+1)}{s^2}v+b_p(s)|v|^p
\big)(1+h,\cdot)
\big\|_\infty\\
&\leq
C\ve h\bigl(\|\na^2f\|_\infty+\|\na g\|_\infty\bigr)
+\frac{h}{2}
\bigl(\lambda(\lambda+1)M_T+M_T^p\bigr)
\longrightarrow0  \qquad \text{as}\, h\to0^+
\end{aligned}
\end{equation*}
which derives $\p_t v(1,x)=\ve g(x)$.
\end{proof}

\section{Proof of Theorem~\ref{thm:main}}\label{Sec5}

Based on Proposition~\ref{prop:iteration},
we now prove Theorem~\ref{thm:main}.

\begin{proof}[Proof of Theorem~\ref{thm:main}]
By $a=\frac{\mu}{2}<\frac{3}{2}$ and $p>\frac{5}{3}$,
we may choose $q$ such that
\begin{equation}\label{qq:156}
\frac{5}{3}<q<\min\big\{p,1+\frac{2}{\mu}\big\},
\end{equation}
and then $\kappa$ can be chosen as in \eqref{qq:67}.

For $f=u_0$ and $g=u_1+au_0$, it follows from
Lemma~\ref{lem:homogeneous-majorant} that there exist
nonnegative radially symmetric functions
$F,G\in C_0^\infty(B(0,1))$ such that
\begin{equation}\label{qq:157}
H(F,G)\geq|H(f,g)|.
\end{equation}
Set
\begin{equation}\label{qq:158}
\phi_0=F,\quad \phi_1=G-aF,\quad \delta=\ve.
\end{equation}
For $0<\ve\leq\delta_*$, Proposition~\ref{prop:radial}
and Lemmas~\ref{lem:compatibility}--\ref{lem:radial-nonnegative}
provide a radial auxiliary solution $W=t^aU\in X_\kappa$
satisfying
\begin{equation}\label{qq:159}
\begin{aligned}
W&=\ve H(F,G)+I\big(
\frac{\lambda(\lambda+1)}{s^2}W+b_q(s)W^q
\big),\\
&\qquad\quad0\leq W(t,x)\leq\frac{C\ve}{1+t}.
\end{aligned}
\end{equation}
Let
\begin{equation}\label{qq:160}
0<\ve\leq\ve_0  =\min\{\delta_*,C^{-1}\}.
\end{equation}
Then, for $t\geq1$ and $x\in\R^3$, one can obtain
\begin{equation}\label{qq:161}
0\leq t^{-a}W(t,x)
\leq C\ve\,t^{-a}(1+t)^{-1}\leq1.
\end{equation}
By \eqref{qq:22}, \eqref{qq:156} and \eqref{qq:161},
it holds that
\begin{equation}\label{qq:162}
\begin{aligned}
b_p(t)W^p
&=t^{-a(p-1)}W^p=t^{-a(q-1)}W^q\bigl(t^{-a}W\bigr)^{p-q}
\leq b_q(t)W^q.
\end{aligned}
\end{equation}
This, together with \eqref{qq:157}, \eqref{qq:159}
 and  \eqref{qq:36}, yields
\begin{equation}\label{qq:163}
W\geq\ve|H(f,g)|
+I\big(
\frac{\lambda(\lambda+1)}{s^2}W+b_p(s)W^p
\big).
\end{equation}
Thus \eqref{qq:121} is satisfied.
By Proposition~\ref{prop:iteration}, there exists a
global solution $v\in C([1,\infty)\times\R^3)$
of \eqref{qq:25} satisfying
\begin{equation}\label{qq:164}
|v(t,x)|\leq W(t,x)
\lesssim \frac{\ve}
{\langle t+|x|\rangle\langle t-|x|\rangle^{\kappa-1}}
\lesssim \frac{\ve}{1+t}.
\end{equation}
We next establish the regularity of $v$.
For any fixed $T>1$ and $t,s\in[1,T]$,
it follows from $W\in X_\kappa$ and \eqref{qq:164} that
\begin{equation}\label{qq:165}
\supp v(t,\cdot)\subset\overline{B(0,t)}
\subset\overline{B(0,T)}.
\end{equation}
By \eqref{qq:165} and the continuity of $v$
on the compact set $[1,T]\times\overline{B(0,T)}$,
one then has
\begin{equation}\label{qq:167}
\|v(t,\cdot)-v(s,\cdot)\|_\infty
=\sup_{|x|\leq T}|v(t,x)-v(s,x)|
\longrightarrow0\quad \text{as}\,\, s\to t
\end{equation}
and
\begin{equation}\label{qq:168}
\|v(t,\cdot)-v(s,\cdot)\|_2
\leq |B(0,T)|^{\frac{1}{2}}
\|v(t,\cdot)-v(s,\cdot)\|_\infty
\longrightarrow0\quad \text{as}\,\, s\to t.
\end{equation}
Hence $v\in C([1,T];L^2(\R^3)\cap L^\infty(\R^3))$.
Moreover, by \eqref{qq:127}, \eqref{qq:164}
and the mean value theorem, it holds that
\begin{equation}\label{qq:169}
\begin{aligned}
\bigl\||v(t)|^p-|v(s)|^p\bigr\|_2
&\leq\big(
\int_{\R^3}p^2M_T^{2(p-1)}
\bigl||v(t,x)|-|v(s,x)|\bigr|^2\dd x
\big)^{\frac{1}{2}}\\
&\leq pM_T^{p-1}\|v(t)-v(s)\|_2.
\end{aligned}
\end{equation}
Similarly,
\[
\bigl\||v(t)|^p\bigr\|_2
\leq M_T^{p-1}\|v(t)\|_2.
\]
Together with \eqref{qq:22}, these estimates yield
\begin{equation}\label{qq:170}
\begin{aligned}
&\big\|
\frac{\lambda(\lambda+1)}{t^2}v(t)+b_p(t)|v(t)|^p
-\frac{\lambda(\lambda+1)}{s^2}v(s)-b_p(s)|v(s)|^p
\big\|_2\\
&\leq
\frac{\lambda(\lambda+1)}{t^2}\|v(t)-v(s)\|_2
+b_p(t)\bigl\||v(t)|^p-|v(s)|^p\bigr\|_2\\
&\quad+
\lambda(\lambda+1)|t^{-2}-s^{-2}|\|v(s)\|_2
+|b_p(t)-b_p(s)|\bigl\||v(s)|^p\bigr\|_2\\
&\leq
\bigl(\lambda(\lambda+1)+pM_T^{p-1}\bigr)
\|v(t)-v(s)\|_2\\
&\quad+
\big[
\lambda(\lambda+1)|t^{-2}-s^{-2}|
+M_T^{p-1}|b_p(t)-b_p(s)|
\big]\|v(s)\|_2
\longrightarrow0\quad \text{as}\,\, s\to t.
\end{aligned}
\end{equation}
Consequently,
\[
\frac{\lambda(\lambda+1)}{t^2}v+b_p(t)|v|^p
\in C([1,T];L^2(\R^3)).
\]
Then by \eqref{qq:142}, we arrive at
\begin{equation}\label{qq:171}
v\in C([1,T];H^1(\R^3))
\cap C^1([1,T];L^2(\R^3)).
\end{equation}

In order to obtain
$v\in C([1,T];H^2(\R^3))$, it remains to prove
$|v|^p\in C([1,T];H^1(\R^3))$.
Choosing
$\rho\in C_0^\infty(B(0,1))$ with
$\rho\geq0$ and $\int_{\R^3}\rho(y)\dd y=1$.
For $\eta>0$ and $t\in[1,T]$, set
\[
v_\eta(t,x)=\int_{\R^3}\rho(y)v(t,x-\eta y)\dd y.
\]
We next prove that as $\eta\to0^+$,
\begin{equation}\label{qq:175}
\begin{gathered}
v_\eta(t)\longrightarrow v(t)
\quad\text{in }H^1(\R^3)\cap L^\infty(\R^3),\\
\|v_\eta(t)\|_\infty
\leq\|v(t)\|_\infty\leq M_T.
\end{gathered}
\end{equation}

Indeed, it follows from \eqref{qq:165} and the continuity
of $v$ that $v(t,\cdot)$ is uniformly continuous on $\R^3$.
Hence, one has
\[
\begin{aligned}
&\|v_\eta(t)-v(t)\|_\infty=\big\|
\int_{\R^3}\rho(y)
\bigl[v(t,\cdot-\eta y)-v(t)\bigr]\dd y
\big\|_\infty\leq\sup_{|y|\leq1}
\|v(t,\cdot-\eta y)-v(t)\|_\infty
\longrightarrow0.
\end{aligned}
\]
Moreover, from \eqref{qq:171}, we may apply
the approximation property of mollifiers to obtain
\[
\|v_\eta(t)-v(t)\|_{H^1}\longrightarrow0.
\]
In addition, we deduce from \eqref{qq:127} and \eqref{qq:164} that
\[
\begin{aligned}
\|v_\eta(t)\|_\infty
&\leq\int_{\R^3}\rho(y)
\|v(t,\cdot-\eta y)\|_\infty\dd y=\|v(t)\|_\infty\int_{\R^3}\rho(y)\dd y=\|v(t)\|_\infty\leq M_T.
\end{aligned}
\]
Thus \eqref{qq:175} holds.
We next prove $|v(t)|^p\in H^1(\R^3)$.
By \eqref{qq:175} and the mean value
theorem, one has that for $\eta\to0^+$,
\begin{equation}\label{qq:177'}
\bigl\||v_\eta(t)|^p-|v(t)|^p\bigr\|_2
\leq pM_T^{p-1}\|v_\eta(t)-v(t)\|_2
\longrightarrow0
\end{equation}
and
\begin{equation}\label{qq:177}
\begin{aligned}
&\bigl\|\na|v_\eta(t)|^p
-p|v(t)|^{p-2}v(t)\na v(t)\bigr\|_2\\
&=p\Bigl\|
|v_\eta(t)|^{p-2}v_\eta(t)
\bigl(\na v_\eta(t)-\na v(t)\bigr)+
\bigl(|v_\eta(t)|^{p-2}v_\eta(t)
-|v(t)|^{p-2}v(t)\bigr)\na v(t)
\Bigr\|_2\\
&\leq pM_T^{p-1}\|\na v_\eta(t)-\na v(t)\|_2+p
\bigl\||v_\eta(t)|^{p-2}v_\eta(t)
-|v(t)|^{p-2}v(t)\bigr\|_\infty
\|\na v(t)\|_2
\longrightarrow0,
\end{aligned}
\end{equation}
where in the last limit we have used the uniform continuity of
$z\mapsto |z|^{p-2}z$ on $[-M_T,M_T]$ for $p>1$.
Based on \eqref{qq:177'} and \eqref{qq:177}, we now derive
$\na|v(t)|^p=p|v(t)|^{p-2}v(t)\na v(t)$
in the weak sense.
For any $\varphi\in C_0^\infty(\R^3)$ and $j=1,2,3$,
one has
\begin{equation}\label{qq:178}
\begin{aligned}
\int_{\R^3}|v(t)|^p\p_j\varphi\dd x
&=\lim_{\eta\to0^+}
\int_{\R^3}|v_\eta(t)|^p\p_j\varphi\dd x\\
&=-\lim_{\eta\to0^+}
\int_{\R^3}
p|v_\eta(t)|^{p-2}v_\eta(t)\p_jv_\eta(t)\varphi\dd x\\
&=-\int_{\R^3}
p|v(t)|^{p-2}v(t)\p_jv(t)\varphi\dd x.
\end{aligned}
\end{equation}
Thus $|v(t)|^p\in H^1(\R^3)$ and
\begin{equation}\label{qq:179}
\begin{gathered}
\na|v(t)|^p=p|v(t)|^{p-2}v(t)\na v(t),\\
\bigl\|\na|v(t)|^p\bigr\|_2
\leq pM_T^{p-1}\|\na v(t)\|_2.
\end{gathered}
\end{equation}

We now prove that $|v|^p\in C([1,T];H^1(\R^3))$.
For $t,s\in[1,T]$, it follows from \eqref{qq:179}, \eqref{qq:167}, \eqref{qq:171} and the uniform continuity of
$z\mapsto |z|^{p-2}z$ on $[-M_T,M_T]$ that
\begin{equation}\label{qq:181}
\begin{aligned}
&\bigl\|\na|v(t)|^p-\na|v(s)|^p\bigr\|_2\\
&=p\bigl\|
|v(t)|^{p-2}v(t)\bigl(\na v(t)-\na v(s)\bigr)
+\bigl(|v(t)|^{p-2}v(t)-|v(s)|^{p-2}v(s)\bigr)
\na v(s)
\bigr\|_2\\
&\leq pM_T^{p-1}\|\na v(t)-\na v(s)\|_2
+p\bigl\||v(t)|^{p-2}v(t)-|v(s)|^{p-2}v(s)\bigr\|_\infty
\|\na v(s)\|_2\\
&\longrightarrow0 \quad\text{as}\,\, s\to t.
\end{aligned}
\end{equation}
Together with \eqref{qq:22}, \eqref{qq:169}
and \eqref{qq:171}, this proves that
\begin{equation}\label{qq:182}
|v|^p,\quad
\frac{\lambda(\lambda+1)}{t^2}v+b_p(t)|v|^p
\in C([1,T];H^1(\R^3)).
\end{equation}
By problem \eqref{qq:25}, we have
\begin{equation}\label{qq:183}
v\in C([1,T];H^2(\R^3))
\cap C^1([1,T];H^1(\R^3)).
\end{equation}
Moreover, by \eqref{qq:25}, \eqref{qq:182} and \eqref{qq:183},
one has
\begin{equation}\label{qq:184}
\p_t^2v
=\Delta v+\frac{\lambda(\lambda+1)}{t^2}v
+b_p(t)|v|^p
\in C([1,T];L^2(\R^3)).
\end{equation}
Consequently,
\begin{equation}\label{qq:186}
v\in C([1,\infty);H^2(\R^3))
\cap C^1([1,\infty);H^1(\R^3))
\cap C^2([1,\infty);L^2(\R^3)).
\end{equation}
By returning to $u=t^{-a}v$ with $\mu=2a$ and
$\lambda=a-1$, we conclude from \eqref{qq:25}
and \eqref{qq:186} that $u$ satisfies \eqref{qq:1}
with the regularity stated in \eqref{qq:5}.
Moreover, for $t\geq1$, it follows from \eqref{qq:164} that
\begin{equation}\label{qq:166}
\|u(t,\cdot)\|_\infty
=t^{-a}\|v(t,\cdot)\|_\infty
\leq C\ve\,t^{-\frac{\mu}{2}}(1+t)^{-1},
\end{equation}
which proves \eqref{qq:6}.
\end{proof}

	\vskip 0.1 true cm
	\noindent
	\textbf{Data availability statement.}
	All mathematical arguments are contained in this article
	and the cited references. No datasets were generated or
	analysed during this study.

\end{document}